\documentclass{article}

\usepackage{amsmath}
  \usepackage{paralist}
  \usepackage{graphics} %% add this and next lines if pictures should be in esp format
  \usepackage{epsfig} %For pictures: screened artwork should be set up with an 85 or 100 line screen
  \usepackage{amssymb,amsmath}
  \usepackage{amsfonts}
  \usepackage{mathrsfs,enumerate}
  \usepackage[colorlinks=true]{hyperref}

\newtheorem{theorem}{Theorem}[section]

\newcommand{\ii}{{\mbox{\boldmath$i$}}}

\newcommand{\uu}{{\mbox{\boldmath$u$}}}
\newcommand{\vv}{{\mbox{\boldmath$v$}}}
\newcommand{\ww}{{\mbox{\boldmath$w$}}}

\newcommand{\tauV}{{\kern-3pt\tau}}

\renewcommand{\AA}{\mathscr{A}}

\newcommand{\XX}{{\mbox{\boldmath$X$}}}

\newcommand{\nnu}{{\mbox{\boldmath$\nu$}}}
\newcommand{\snnu}{{\mbox{\scriptsize\boldmath$\nu$}}}

\newcommand{\eeta}{{\mbox{\boldmath$\eta$}}}

\newcommand{\R}{\mathbb{R}}

\newcommand{\Z}{\mathbb{Z}}
\newcommand{\T}{\mathbb{T}}

\newcommand{\Leb}[1]{{\mathscr L}^{#1}}      % Misura di Lebesgue

\newcommand{\Probabilities}[1]{\mathscr P(#1)}          % misure di probabilita'
\newcommand{\BorelSets}[1]{\mathcal B(#1)}

\newcommand{\SDiff}[1]{{\rm SDiff}(#1)}
\renewcommand{\div}{{\rm div}}

\title{Lecture notes on
variational models\\ for incompressible Euler equations}

\author{Luigi Ambrosio\footnote{Scuola Normale Superiore di Pisa, Piazza dei Cavalieri 7,
56126 Pisa, Italy
(\texttt{l.ambrosio@sns.it})}\ and Alessio
Figalli\footnote{Department of Mathematics, The University of Texas
at Austin, Austin TX 78712, USA (\texttt{figalli@math.utexas.edu})}}

\begin{document}
\maketitle

\begin{abstract}
These notes briefly summarize the lectures for the Summer School "Optimal transportation: Theory and applications" held by the second author in Grenoble during the week of June 22-26, 2009.
Their goal is to describe some recent results on
 Brenier's variational models for incompressible Euler equations \cite{amfiga1,amfiga,berfigsan}.
\end{abstract}

\section{Euler incompressible equations and Arnold geodesics}

Let $D$ denote either a bounded domain of $\R^d$ or the $d$-dimentional torus $\T^d:=\R^d/\Z^d$.
We consider an incompressible fluid moving inside
$D$ with velocity $\uu$. The Euler equations
for $\uu$ describe the evolution in time of the velocity field, and are
given by
$$
\left\{
\begin{array}{ll}
\partial_t \uu + (\uu \cdot \nabla)\uu= - \nabla p & \text{in }[0,T] \times D,\\
\div\,\uu=0 & \text{in }[0,T] \times D,\\
\end{array}
\right.
$$
coupled with the boundary condition
$$
\uu \cdot \nu=0 \qquad \text{on }[0,T] \times \partial D
$$
when $D \neq \T^d$.
Here $p$ is the pressure field, and arises as a Lagrange multiplier for the
divergence-free constraint on the velocity
$\uu$.

If $\uu$ is smooth we can write the above equations in Lagrangian
coordinates: let $g$ denote the flow map of $\uu$, that is
$$
\left\{
\begin{array}{l}
\dot g(t,a)=\uu(t,g(t,a)),\\
g(0,a)=a.
\end{array}
\right.
$$
By the incompressibility condition, and the classical differential
identity
$$
\frac{d}{dt}{\rm det}\nabla_a g(t,a)=\div\,\uu(t,g(t,a))\,{\rm
det}\nabla_a g(t,a),
$$
(here and in the sequel $\div$ denotes the spatial divergence of
a possibly time-dependent vector field)
we get ${\rm det}\nabla_a g(t,a)\equiv 1$. This means that
$g(t,\cdot):D \to D$ is a measure-preserving diffeomorphism of $D$:
$$
g(t,\cdot)_\#\mu_D=\mu_D \qquad \bigl(\text{i.e.
$\mu_D(g(t,\cdot)^{-1}(E))=\mu_D(E)\,\,\,\forall E\bigr)$.}
$$
Here and in the sequel $f_\#\mu$ is the push-forward of a Borel
measure $\mu$ through a map $f:X\to Y$ (i.e. $\int_Y
\phi\,df_\#\mu=\int_X \phi\circ f\,d\mu$ for all Borel bounded functions $\phi:Y\to\R$), and $\mu_D$ is the volume
measure of $D$, renormalized by a constant so that $\mu_D(D)=1$.

Writing Euler's equations in terms of $g$ we obtain an ODE for $t \mapsto g(t)$ in the
space $\SDiff{D}$ of measure-preserving smooth diffeomorphisms of
$D$:
\begin{equation}\label{ODiff} \left\{
\begin{array}{ll}
\ddot g(t,a) =-\nabla p\left(t,g(t,a)\right) & \text{$(t,a)\in [0,T] \times D$,}\\
g(0,a)=a &\text{$a\in D$,}\\
g(t,\cdot)\in \SDiff{D} &\text{$t\in [0,T]$.}
\end{array}
\right.
\end{equation}

\subsection{Weak solutions to Euler's equations}

In the case $d=2$, existence of distributional solutions can be proved through the
\emph{vorticity} formulation: setting $\omega_t(\cdot)={\rm
curl\,}\uu(t,\cdot)$, so that
$\uu(t,\cdot)=\nabla^\perp\Delta^{-1}\omega_t$, the Euler equations
can be read as follows:
$$
\frac{d}{dt} \omega_t(x)+{\rm div\,}\bigl(\omega_t(x)\uu(t,x)\bigr)=0.
$$
Formally, this equation preserves all $L^p$ norms of solutions, and
indeed existence is not hard to obtain if $\omega_0\in L^p$ for $1<
p\leq\infty$. Delort improved the existence theory up to $L^1$ or
measure initial conditions $\omega_0$ whose positive (or negative)
part is absolutely continuous, and it is still open the problem of
getting a solution for all measure initial data. As shown by Yudovitch \cite{yudo1,yudo2}
uniqueness holds for $p=\infty$, while it is still open in all the other
cases.

In the case $d>2$ much less is known: \emph{no} general global existence
results of distributional solutions is presently available.

\subsection{Arnold's geodesic interpretation}

At least formally, one can view the space $\SDiff{D}$ of measure-preserving
diffeomorphisms of $D$ as an infinite-dimensional manifold with the
metric inherited from the embedding in $L^2(D;\R^d)$, and with tangent
space made by the divergence-free vector fields.  Using this
viewpoint, Arnold interpreted the ODE (\ref{ODiff}), and therefore
Euler's equations, as a \emph{geodesic} equation on $\SDiff{D}$.
Therefore one can look for solutions of Euler's equations on $[0,1]\times D$ by
minimizing the Action functional
$$
\AA(g):=\int_0^1\int_D \frac{1}{2} |\dot g(t,x)|^2\,d\mu_D(x)\,dt
$$
among all paths $g(t,\cdot):[0,1]\to\SDiff{D}$ with $g(0,\cdot)=f$
and $g(1,\cdot)=h$ prescribed (typically, by right invariance, $f$
is taken as the identity map $\ii$). Ebin and Marsden proved in
\cite{ebinmarsden} that this problem has indeed a unique solution
when $h\circ f^{-1}$ is sufficiently close, in a strong Sobolev
norm, to $\ii$. We shall denote by $\delta(f,h)$ the Arnold distance
in $\SDiff{D}$ induced by this minimization problem.

Of course, this variational problem differs from Euler's problem,
because the initial and final diffeomorphisms, and not the initial
velocity, are prescribed. Nevertheless, the investigation of this
problem leads to difficult and still not completely understood
questions (typical of Calculus of Variations) namely:

\begin{itemize}
\item[(a)] Necessary and sufficient optimality conditions;
\item[(b)] Regularity of the pressure field;
\item[(c)] Regularity of (relaxed) curves with minimal length.
\end{itemize}

Before describing some of the main contributions in this field, let us recall
some ``negative'' results that motivate somehow the necessity of
relaxed formulations of this minimization problem.

\subsection{Non-attainment and non-existence results}

Shnirelman \cite{shnir1,shnir2} found the example of a map $\bar g\in\SDiff{[0,1]^2}$ which cannot
be connected to $\ii$ by a path with finite action, i.e.
$\delta(\ii,\bar g)=+\infty$.
Furthermore, he proved that for
$h\in\SDiff{[0,1]^3}$ of the form
$$
h(x_1,x_2,x_3)=(\bar g_1(x_1,x_2),\bar g_2(x_1,x_2),x_3),\quad\text{with}
\quad (\bar g_1,\bar g_2)=\bar g\text{ as above,}
$$
$\delta(\ii,h)$ is not attained, i.e. no minimizing path
between $\ii$ and $h$ exists (although there exist paths with a
finite action).
This fact can be easily explained as follows (see also \cite[Paragraph 1.3]{brenier4}): since there is no two dimensional path
with finite action connecting $\ii$ to $\bar g$ while in $3$ dimension it is known that
the minimal action is finite \cite{shnir1}, if a minimizing path $t \mapsto g(t)$ exists
then it has a non-trivial third component, i.e. $g_3(t,x)\not\equiv x_3$.
Set $\eta(x_3):=\min\{2x_3,2-2x_3\}$, and let
$\uu$ denote the velocity field associated to $g$, i.e. $\uu=\dot g\circ g^{-1}$.
Then it is easily seen that the velocity field
$$
\tilde \uu(x_1,x_2,x_3):=\left\{\begin{array}{l}
\uu_1(x_1,x_2,\eta(x_3))\\
\uu_2(x_1,x_2,\eta(x_3))\\
\frac{1}{2}\uu_3(x_1,x_2,\eta(x_3))
\end{array}
\right.
$$
induces a path $\tilde g$ which still joins $\ii$ to $h$, but with strictly less action (since $\uu_3$
is not identically zero).
This contradicts the minimality of $g$, and proves that there is no minimizing path between $\ii$ and $h$. (See also \cite[Paragraph 1.3]{brenier4}.)

Let us point out that the above argument shows that minimizing
sequences exhibit oscillations on small
scales, and strongly suggest the analysis of weak solutions.

\subsection{Time discretization, minimal projection and optimal transport}
Before describing the concept of relaxed solutions to the Euler
equations introduced by Brenier, let us first see what happens when
one tries to attack the above variational problem by
time-discretization: assume $D\subset \R^d$, and fix $g_0,g_1 \in
\SDiff{D}$. We want to find the ``midpoint" $g_{1/2}$ between $g_0$
and $g_1$, that is we consider
$$
\min_{g \in \SDiff{D}} \Bigl\{\frac12 \|g-g_0\|_{L^2(D;\R^d)}^2 +\frac12\|g_1-g\|_{L^2(D;\R^d)}^2 \Bigr\}.
$$
Up to rearranging the terms and removing all the quantities independent on $g$,
the above problem is equivalent to minimize
$$
\min_{g \in \SDiff{D}}\Bigl\|g-\frac{g_0+g_1}{2}\Bigr\|_{L^2(D;\R^d)}^2,
$$
i.e. we have to find the $L^2$-projection on $\SDiff{D}$ of the function $\frac{g_0+g_1}{2}\in L^2(D;\R^d)$.
Since the set $\SDiff{D}$ is neither closed nor convex, no classical theory is available to ensure the
existence of such projection.

In order to make the problem more treatable, let us close $\SDiff{D}$: as shown for instance
in \cite{brengang}, if $D=[0,1]^d$ or $D=\T^d$ then the $L^2$-closure of $\SDiff{D}$ in $L^2(D;\R^d)$
coincides with the space $S(D)$ of measure preserving maps:
$$
S(D):=\left\{g:D\to D:\
\mu_D(g^{-1}(A))=\mu_D(A)\,\,\forall A\in\BorelSets{D}\right\}.
$$
Then the general problem we want to study becomes the following:
given $h\in L^2(D;\R^d)$, solve
\begin{equation}
\label{eq:min s}
\min_{s \in S(D)} \int_D |h-s|^2\,d\mu_D.
\end{equation}
As in the classical optimal transport problem, one can consider the following Kantorovich
relaxation: denoting by $\Pi(\R^d)$ the set of probability measures on $\R^d\times\R^d$
with first marginal $\mu_D$ and second marginal $\nu:=h_\#\mu_D$, we minimize
\begin{equation}
\label{eq:min gamma}
\min_{\gamma \in \Pi(\R^d)} \int_{\R^d\times\R^d} |x-y|^2\,d\gamma(x,y).
\end{equation}
Assume the \textit{non-degeneracy condition} $\nu\ll
dx$. Then we can apply the classical theory of optimal transport
with quadratic cost for the problem of sending $\nu$ onto $\mu_D$
\cite{breniertransp}: there exists a unique optimal transport map
$\nabla \phi:\R^d\to\R^d$ such that $(\nabla\phi)_\#\nu=\mu_D$.
Moreover the unique optimal measure $\bar\gamma$ which solves
\eqref{eq:min gamma} is given by
$$
\bar \gamma=(\nabla\phi \times {\rm Id})_\# \nu.
$$
Then it is easily seen that the map
$$
\bar s:= \nabla \phi \circ h
$$
belongs to $S(D)$ and uniquely solves \eqref{eq:min s}
(see \cite{breniertransp} or \cite[Chapter 3]{villani} for more details).

\section{Relaxed solutions}
In the last paragraph we have seen how the attempt of attacking
Arnold's geodesics problem by time discretization leads to study the
existence of the $L^2$-projection onto $\SDiff{D}$, and that the
projection of a function $h$ onto its closure $S(D)$ exists and is
unique whenever $h$ satisfies a non-degeneracy condition. Instead of
going on with this strategy, we now want to change point of view,
attacking the problem by a relaxation in ``space".

Two levels of relaxation can be imagined: the first one is to relax
the smoothness and injectivity constraints, and this leads to the
definition of the space $S(D)$ of \emph{measure-preserving maps}.
However, we will see that a second level is necessary, giving up the idea
that $g(t,\cdot)$ is a map, but allowing it to be a \emph{measure
preserving plan} (roughly speaking, a multivalued map). This leads
to the space
$$
\Gamma(D):=\left\{
\eta\in\Probabilities{D\times D}:\
\eta(A\times D)=\mu_D(A)=\eta(D\times A)\,\,\,\forall A\in\BorelSets{D}
\right\}.
$$
The space $S(D)$ ``embeds'' into $\Gamma(D)$ considering
$$
S(D)\ni g\mapsto (\ii\times g)_\#\mu_D\in\Gamma(D).
$$
Conversely, any $\eta\in\Gamma(D)$ concentrated on a graph is
induced by a map $g\in S(D)$.

Even from the Lagrangian viewpoint, it is natural to follow the path
of each particle, and to relax the smoothness and injectivity
constraints, allowing fluid paths to split, forward or backward in
time. These remarks led in 1989 Brenier to the following model \cite{brenier1}: let
$$
\Omega(D):=C\left([0,1];D\right),
\qquad e_t(\omega):=\omega(t),\,\,t\in [0,1].
$$
Then, denoting by $\Probabilities{\Omega(D)}$ the family of
probability measures in $\Omega(D)$, we minimize the action
functional
$$
\AA(\eeta):=\int_{\Omega(D)}\frac{1}{2}\int_0^1|\dot\omega|^2\,dt\,d\eeta(\omega),
\qquad\eeta\in\Probabilities{\Omega(D)}
$$
with the endpoint and incompressibility constraints
$$
(e_0,e_1)_\#\eeta=(\ii\times h)_\#\mu_D,
\qquad (e_t)_\#\eeta=\mu_D\,\,\,\forall t\in [0,T].
$$
In Brenier's model, a flow is modelled by a random path with some
constraints on the expectations of this path.
As we will see below, this problem can be
recast in the optimal transportation framework, dealing properly
with the incompressibility constraint.

Classical flows $g(t,a)$ induce generalized ones, with the same
kinetic action, via the relation $\eeta=(\Phi_g)_\#\mu_D$, with
$$ \Phi_g:D\to\Omega(D),\qquad\quad\Phi_g(a):=g(\cdot,a).$$

In this relaxed model, some obstructions of the original one
disappear: for instance, in the case $D=[0,1]^d$ or $D=\T^d$
it is always possible to connect any couple of measure preserving diffeomorphism by a path
with action less than $\sqrt{d}$. Actually, this allows to prove that finite-action paths exist
in many situation: as shown in \cite[Theorem 3.3]{amfiga1},
given a domain $D$ for which there exists a bi-Lipschitz measure-preserving diffeomorphism $\Phi:D\to [0,1]^d$, 
by considering composition of generalized flows with $\Phi$ one can easily constructs a generalized flow with finite action
between any $h_0,h_1\in\SDiff{D}$. 
 Moreover, standard
compactness/lower semicontinuity arguments in the space
$\Probabilities{\Omega(D)}$ provide existence of generalized flows
with minimal action.

\subsection{Eulerian-Lagrangian model}
Coming back to the relaxed model described above,
we observe that the endpoint constraint
$(e_0,e_1)_\#\eeta=(\ii\times h)_\#\mu_D$ cannot be modified to deal
with the more general problem of connecting $f\in S(D)$ to $h\in
S(D)$: indeed, by right invariance, this is clear only if $f$ is
invertible (in this case, one looks for the optimal connection
between $\ii$ and $h\circ f^{-1}$). These remarks led to a more
general model, which allows to connect $\eta=\eta_a\otimes\mu_D$ to
$\gamma=\gamma_a\otimes\mu_D$ \cite{amfiga}. (Here we are disintegrating both the
initial and final plan with respect to the first variable.) The idea, which appears first in
Brenier's \emph{Eulerian-Lagrangian} model \cite{brenier4} is to
``double'' the state space, adding to the Eulerian state space $D$ a
Lagrangian state space $A$. Even though $A$ could be thought as an
identical copy of $D$, it is convenient to denote it by a different
symbol.

Let
$$
\Omega^*(D):=\Omega(D)\times A.
$$
Then, consider probability measures $\eeta=\eeta_a\otimes\mu_D$
in $\Omega^*(D)$: this means that $\eeta$ has $\mu_D$ as second
marginal, and that
$$
\int
\phi(\omega,a)\,d\eeta(\omega,a)=\int_A\biggl(\int_{\Omega(D)}\phi(\omega,a)\,
d\eeta_a\biggr)\,d\mu_D(a)
$$
for all bounded Borel functions $\phi$ on $\Omega^*(D)$.

Again, one minimizes the action
$$
\AA(\eeta):=\int_{\Omega^*(D)}\frac{1}{2}\int_0^1|\dot\omega|^2\,dt\,d\eeta(\omega,a)
$$
with the incompressibility constraint $(e_t)_\#\eeta=\mu_D$ for all
$t$ (here $e_t(\omega,a)=\omega(t)$) and
the family of endpoint constraints:
$$
(e_0)_\#\eeta_a=\gamma_a,\quad
(e_1)_\#\eeta_a=\eta_a\qquad\text{for $\mu_D$-a.e $a\in D$.}
$$
As in the previous section, we are using
$\eta_a\otimes\mu_D$ and $\gamma_a\otimes\mu_D$ to denote the disintegrations
of $\eta$ and $\gamma$ respectively.

Denoting by $\overline\delta(\eta,\gamma)^2$ the minimal action, it
turns out that one can define natural operations of
\emph{reparameterization}, \emph{restriction} and
\emph{concatenation} in this class of flows. These imply
that $(\overline\delta,\Gamma(D))$ is a metric space.

Indeed, it is proved in \cite{amfiga1} that it is \emph{complete} and a
\emph{length} space, whose convergence is stronger than weak
convergence in $\Probabilities{D\times D}$.

\subsection{Motivation for the extension to $\Gamma(D)$}\label{sbrenier}

Even for deterministic initial and final data, there exist examples
of minimizing geodesics $\eeta$ that are \emph{not deterministic} in
between: this means that $(e_0,e_t)_\#\eeta\in\Gamma(D)\setminus
S(D)$, $t\in (0,1)$.

To show this phenomenon, consider the problem of connectingup to additive
constants in $D=B_1(0)\subset\R^2$ the identity map $\ii$ to
$-\ii$. For convenience, up to a reparameterization, we can
choose the time interval as $[0,\pi]$. Two classical solutions are
$$
[0,\pi]\ni t\mapsto (x_1\cos \pm t+x_2\sin \pm t,x_1\sin\pm
t+x_2\cos\pm t),
$$
corresponding to a clockwise and an anti-clockwise rotation.

On the other hand, one can consider the family of maps
$\omega_{x,\theta}$ connecting $x$ to $-x$
\begin{equation}\label{paths}
\omega_{x,\theta}(t):=x\cos t+\sqrt{1-|x|^2}(\cos\theta,\sin\theta)\sin t
\qquad\theta\in (0,\pi)
\end{equation}
and define $\eeta:={(\omega_{x,\theta})_\sharp}
\bigl(\frac{1}{2\pi^2}\Leb{2}\lfloor D\times\Leb{1}\lfloor
(0,2\pi)\bigr)$.

It turns out that $\eeta$ is optimal as well, and non-deterministic
in between.
Moreover, as shown in \cite{berfigsan}, it is possible to construct infinitely many other
solutions to the above minimization problem which are not induced by maps.
For instance, one can split the measure $\eeta$ above as $\frac{1}{2}\bigl(\eeta_+ + \eeta_-\bigr)$,
where $\eeta_+$ consists of the curves such that
$(\cos\theta,\sin\theta)\cdot x^\perp \geq 0,$
and $\eeta_-$ consists of the curves such that
$(\cos\theta,\sin\theta)\cdot x^\perp \leq 0,$
where $x^\perp=(x_2,-x_1)$,
and the two flows $\eeta_+$ and $\eeta_-$
can be shown to be still incompressible (see \cite[Paragraph 4.1]{berfigsan}).
We will say more about these important examples later on,
as more results on the theory will be available.

\section{The pressure field}

Brenier proved in \cite{brenier2} a surprising result: even though
geodesics are not unique in general, given the initial and final
conditions, there is a \emph{unique}, up to an additive
time-dependent constant, pressure field. The pressure field arises
if one relaxes the incompressibility constraint, considering
\emph{almost incompressible} flows $\nnu$. Denoting by
$\rho^{\snnu}$ the density produced by the flow, defined by
$$
(e_t)_\#\nnu=\rho^{\snnu}(t,\cdot)\mu_D \qquad \left( \text{i.e.
$\int
\phi(\omega(t))\,d\eeta(\omega)=\int_D\phi\rho^{\snnu}(t,\cdot)\,d\mu_D$
for all $\phi$} \right),
$$
we say that $\nnu$ is almost incompressible if
$\Vert\rho^{\snnu}-1\Vert_{C^1}\leq 1/2$.

\begin{theorem}[Pressure as a Lagrange multiplier,
\cite{brenier2,amfiga1}]\label{tcol} Let $\eeta$ be optimal between $\eta$
and $\gamma$. There exists a distribution $p\in (C^1)^*$ such that
\begin{equation}\label{dualbrenier}
\AA(\nnu)+\langle p,\rho^{\snnu}-1\rangle \geq\AA(\eeta)
\end{equation}
for all almost incompressible flows $\nnu$ between $\eta$ and
$\gamma$ with $\rho^{\snnu}(t,\cdot)=1$ for $t$ sufficiently close
to $0$ and to $1$.
\end{theorem}

Using this result one can make first variations as follows: given
a smooth field $\ww(t,x)$, vanishing for $t$ close to $0$ and $1$, one can consider
the family $(\XX^t)$ of flow maps
$$
\frac{d}{d\varepsilon}
\XX^t(\varepsilon,x)=\ww(t,\XX^t(\varepsilon,x)),\qquad\XX^t(0,x)=x
$$
and perturb (smoothly) the paths $\omega$ by
$\omega(t)\mapsto\XX^t(\varepsilon,\omega(t))\sim \omega(t)+\varepsilon \,\ww(t,\omega(t))$.
Denoting by
$$
\Phi_\varepsilon:\Omega^*(D)\to\Omega^*(D),\qquad
\Phi_\varepsilon(\omega,a)(t):=\bigl(\XX^t(\varepsilon,\omega(t)),a\bigr),
$$
the induced perturbations in $\Omega^*(D)$, these in
turn induce perturbations
$\eeta_\varepsilon:=(\Phi_\varepsilon)_\#\eeta$ of $\eeta$ which
are almost incompressible. Then, the first variation gives
$$
\int_{\Omega^*(D)}\int_0^1\dot\omega(t)\cdot\frac{d}{dt}
\ww(t,\omega(t))\,dt\,d\eeta(\omega,a)+ \langle
p,\div\,\ww\rangle=0.
$$
This equation uniquely determines $\nabla p$ as a distribution, independently of the chosen
minimizer $\eeta$: indeed, $\eeta$ enters in \eqref{dualbrenier}
only through $\AA(\eeta)$, which obviously is independent of the
chosen minimizer, and so the above equation holds true for \emph{every}
minimizer $\eeta$. Since $\ww$ is arbitrary, the first variation also
leads to a weak formulation of Euler's equations
$$
\partial_t\overline{\vv}_t(x)
+\div\left(\overline{\vv\otimes\vv}_t(x)\right)+\nabla_x p(t,x)=0,
$$
where $\overline\vv_t$ and $\overline{\vv\otimes\vv}_t$ are
implicitly defined by
$$\overline\vv_t\mu_D=(e_t)_\#(\dot\omega(t)\eeta),\qquad
\overline{\vv\otimes\vv}_t\mu_D=
(e_t)_\#(\dot\omega(t)\otimes\dot\omega(t)\eeta).
$$
Observe that in general
$\overline{\vv\otimes\vv}_t
\neq\overline{\vv}_t\otimes\overline{\vv}_t$. Indeed, since these
models allow the passage of many fluid paths at the same point at
the same time (i.e. branching and multiple velocities are possible),
the product $\overline\vv_t(x)\otimes\overline\vv_t(x)$ of the mean
velocity $\overline\vv_t(x)$ with itself might be quite different
from the mean value $\overline{\vv\otimes\vv}_t(x)$ of the product.
This gap precisely marks the difference between genuine
distributional solutions to Euler's equation and ``generalized''
ones (see also \cite[Section 2 and Paragraph 4.4]{berfigsan} for more comments
on this fact).

\section{Necessary and sufficient optimality conditions}

In this section we study necessary and sufficient optimality
conditions for Brenier's variational problem and its extensions.

The basic remark is that \emph{any} Borel integrable function
$q:[0,1]\times D\to\R$ with $\int_D q(t,\cdot)\,d\mu_D=0$ for every
$t\in [0,1]$ induces a null-lagrangian for the minimization problem,
with the incompressibility constraint: indeed
$$
\int_{\Omega^*(D)}\int_0^1 q(t,\omega(t))\,dt\,d\eeta(\omega,a)=
\int_0^1\int_D q(t,x)\,d\mu_D(x)\,dt=0
$$
for any generalized incompressible flow $\eeta$. If we denote by
$$
c_q^{0,1}(x,y):=\inf\left\{\int_0^1\frac{1}{2}|\dot\omega(t)|^2-q(t,\omega(t))\,dt
:\ \omega(0)=x,\,\,\omega(1)=y\right\}
$$
the value function for the Lagrangian ${\mathcal
L}_q(\omega):=\int\frac{1}{2}|\dot\omega(t)|^2-q(t,\omega(t))\,dt$,
we also have
$$
\int_{\Omega^*(D)}\int_0^1\frac{1}{2}|\dot\omega(t)|^2-q(t,\omega(t))\,dt\,d\eeta(\omega,a)
\geq \int_Dc_q^{0,1}(a,h(a))\,d\mu_D(a)
$$
for any incompressible flow $\eeta$ between $\ii$ and $h$. Moreover,
equality holds if and only if $\eeta$-almost every $(\omega,a)$ is a
$c_q^{0,1}$-minimizing path.

The following result, proved in \cite[Section 3.6]{brenier4}, shows that this
lower bound is sharp with $q=p$, if $p$ is sufficiently smooth.

\begin{theorem}
Let $\uu$ be a $C^1$ solution to the Euler equations in
$[0,T]\times D$, whose pressure field $p$ satisfies
$$
T^2\sup_{t\in [0,T]}\sup_{x\in D,\,|\xi|\leq
1}\langle\nabla_x^2p(t,x)\xi,\xi\rangle\leq\pi^2.\leqno (\ast)
$$
Then the measure $\eeta$ induced by $\uu$ via the flow map is
optimal on $[0,T]$.
\end{theorem}

This follows by the fact that the integral paths of $\uu$ satisfy
$\ddot\omega(t)=-\nabla p(t,\omega(t))$, and ($\ast$) implies that
stationary paths for the action are also minimal for ${\mathcal
L}_p$. (This is a consequence of the one dimensional Poincar\'e inequality
$\int_0^T |\dot u(t)|^2\,dt \geq \frac{\pi^2}{T^2}\int_0^T
|u(t)|^2dt$ for all $u:[0,T]\to \R$ such that $\int_0^T u\,dt =0$,
see \cite[Section 5]{brenier1} or \cite[Proposition 3.2]{brenier4} for more details.)

The question investigated in \cite{amfiga1} is: how far are these
conditions from being \emph{necessary}? $C^1$ regularity or even
one sided bounds on $\nabla^2 p$ are not realistic, so one has to look for
necessary (and sufficient) conditions under much weaker regularity
assumptions on
$p$. 

From now on, we restrict for simplicity to the case $D=\T^d$. The following
regularity result for the pressure field has been obtained in
\cite{amfiga}, improving the regularity $\nabla p\in {\mathcal
M}_{\rm loc}\left((0,1)\times\T^d\right)$ obtained in
\cite{brenier4}.

\begin{theorem}
For any $\gamma,\,\eta\in\Gamma(\T^d)$  the unique pressure field given
by Theorem~\ref{tcol} belongs to $L^2_{\rm
loc}\left((0,1);BV(\T^d)\right)$.
\end{theorem}

The above result says in particular that $p$ is a function,
and not just  a distribution. This allows to define the value of $p$
pointwise, which as we will see below will play a key role.

In order to guess the right optimality conditions, we recall that
the two main degrees of freedom in optimal transport problems are:

\smallskip
\noindent
 $\bullet$ In moving mass from $x$ to $y$, the path, or the
family of paths, that should be followed;

\noindent
$\bullet$ The amount of mass that should be moved, on each
such path, from $x$ to $y$.
\smallskip

The second degree of freedom is even more important in situations
when more than one optimal path between $x$ and $y$ is available.
As we will see, both things will depend on ${\mathcal L}_p$. But, since
$p$ is defined only up to negligible sets, the value of the
Lagrangian ${\mathcal L}_p$ on a path $\omega$ is \emph{not}
invariant in the Lebesgue equivalence class; furthermore, no local
pointwise bounds on $p$ are available (remember that $p(t,\cdot)$ is
only a $BV$ function, with $BV$ norm in $L^2_{\rm loc}(0,1)$).
Therefore, as done in \cite{amfiga1}, one has to:

\smallskip
\noindent
 $\bullet$ Define a \emph{precise} representative $\bar p$
in the Lebesgue equivalence class of $p$; it turns out that the
``correct" definition is
$$
\bar p(t,x):=\liminf_{\varepsilon\downarrow
0}p(t,\cdot)\ast\phi_\varepsilon(x),
$$
where $p(t,\cdot)\ast\phi_\varepsilon$ are suitable mollifications
of $p(t,\cdot)$. Of course this definition depends on the choice
of the mollifiers, but we prove that a suitable choice of them
provides a well-behaved (in the sense stated in Theorem~\ref{tollo2}
below) function $\bar p$.

\noindent
$\bullet$ Consider, in the minimization problem, only
paths $\omega$ satisfying
$$
Mp(t,\omega(t))\in L^1_{\rm loc}(0,1),
$$
where $Mp(t,\cdot)$ is a suitable maximal function of $p(t,\cdot)$
(see \cite{amfiga1} for a more precise definition of the maximal
operator).

\smallskip
With these constraints one can talk of \emph{locally minimizing
path} $\omega$ for the Lagrangian ${\mathcal L}_{\bar p}$ and,
correspondingly, define a family of value functions
$$
c^{s,t}_{\bar p}:\T^d\times\T^d\to [-\infty,+\infty],\qquad
[s,t]\subset (0,1),
$$
representing the cost of the minimal connection between $x$ and $y$
in the time interval $(s,t)$:
\begin{multline*}
c^{s,t}_{\bar p}(x,y):=
\inf\biggl\{\int_s^t\frac{1}{2}|\dot\omega(\tau)|^2-\bar
p(\tau,\omega)\,d\tau :\\  \omega(s)=x,\,\,\omega(t)=y,\,\,
Mp(\tau,\omega(\tau))\in L^1(s,t)\biggr\}.
\end{multline*}

With this notation, the following result proved in \cite{amfiga1}
provides necessary and sufficient optimality conditions.

\begin{theorem}\label{tollo2}
Let $\eeta=\eeta_a\otimes\mu_\T$ be an optimal incompressible flow between
$\eta=\eta_a\otimes\mu_\T$ and $\gamma=\gamma_a\otimes\mu_\T$. Then
\begin{itemize}
\item[(i)] $\eeta$ is concentrated on locally minimizing paths for ${\mathcal L}_{\bar p}$;
\item[(ii)] for all intervals $[s,t]\subset (0,T)$, for $\mu_\T$-a.e. $a$, the plan
$(e_s,e_t)_\#\eeta_a$ is $c^{s,t}_{\bar p}$-optimal, i.e.
$$
\int_{\T^d\times\T^d} c^{s,t}_{\bar p}(x,y)\,d(e_s,e_t)_\#\eeta_a\leq
\int_{\T^d\times\T^d} c^{s,t}_{\bar p}(x,y)\,d\lambda
$$
for any $\lambda\in\Probabilities{\T^d\times\T^d}$ having the same marginals
of $(e_s,e_t)_\#\eeta_a$.
\end{itemize}
Conversely, if (i), (ii) hold with $\bar p$ replaced by some function $q$
satisfying
$Mq\in L^1_{\rm loc}\left((0,1);L^1(\T^d)\right)$,
then $\eeta$ is optimal, and $q$ is the pressure field.
\end{theorem}

Notice that an optimal transport problem is trivial if either the
initial or the final measure is a Dirac mass; therefore the second
condition becomes meaningful when either $(e_s)_\#\eeta_a$ or
$(e_t)_\#\eeta_a$ is not a Dirac mass: it corresponds to the case
when $(e_s,\pi_a)_\#\eeta$ is \emph{not} induced by a map, a
phenomenon that cannot be
ruled out, as we discussed in Section~\ref{sbrenier}.
In the example presented in Section~\ref{sbrenier} the
pressure field $p(x)=|x|^2/2$ is smooth and time-independent, but
the initial and final conditions are chosen in such a way that a
continuum of action-minimizing paths (\ref{paths}) between $x$ and
$-x$ exists. As shown in \cite{berfigsan},
there are infinitely many incompressible flows
connecting the identity map $\ii$ to $-\ii$,
which moreover induce infinitely many distributional solutions
to the Euler equations \cite[Paragraph 4.4]{berfigsan}.

The results in \cite{amfiga1} show a connection with the theory of
action-minimizing measures, though
in this case the Lagrangian
$\int_0^1\frac{1}{2}|\dot\omega(t)|^2-\bar p(t,\omega(t))\,dt$ is
possibly non-smooth and not given a priori, but \emph{generated} by
the variational problem itself.

Here we see a nice variation on a classical theme of Calculus of
Variations: a field of (smooth, nonintersecting) \emph{extremals}
gives rise both to \emph{minimizers} and to an incompressible flow
\emph{in phase space}. Here, instead, we have a field of (possibly
nonsmooth, or intersecting) \emph{minimizers} which has to produce
an incompressible flow \emph{in the state space}. This structure
seems to be rigid, and it might lead to new regularity results for
the pressure field.

Let us also recall that, as recently shown in \cite{figmand}, under
a $W^{1,p}$-regularity of the pressure $p$ one can show that
$\eeta$-a.e. $\omega$ solves the Euler-Lagrange equations and
belongs to $W^{2,p}([0,1])\subset C^1([0,1])$. This result is a
first step towards the $BV$ case, where one can still
expect that the minimality of $\eeta$ may allow to prove higher
regularity on the minimizing curves $\omega$ (like $\dot \omega \in
BV$).

%\bibitem{A11}
%     \newblock FirstName MiddleInitial. LastName,
% first name middle initial. and then last name.  Only the first character
% in the paper title is capitalized.
%     \newblock \emph{Title of the paper},
%     \newblock Name of the journal, \textbf{Volume} (Year),
%    StaringPage--EndingPage.
%
% Example of paper with MR number:
%\bibitem{A22} (2082924)
%     \newblock C. Wolf,
%     \newblock \emph{A mathematical model for the propagation of a
%      hantavirus in structured populations},
%     \newblock Discrete Continuous Dynam. Systems - B,
% \textbf{4} (2004), 1065--1089.

 \end{document}